\documentclass{amsart}
\usepackage{eufrak}
\usepackage{amsfonts}
\usepackage{float}
\usepackage{amsmath,amssymb}
\usepackage{lscape}
\usepackage[dvipdfmx]{graphicx, color}
\usepackage[all]{xy}

\newtheorem{thm}{Theorem}[section]
\newtheorem{cor}[thm]{Corollary} 
\newtheorem{lem}[thm]{Lemma} 
\newtheorem{prop}[thm]{Proposition} 
\newtheorem{rem}[thm]{Remark}

\newtheorem{defin}[thm]{Definition}

\newcommand{\A}{\mathcal{A}}
\newcommand{\K}{\mathcal{K}}

\newcommand{\E}{\mathcal{E}}

\newcommand{\R}{\mathbb{R}}

\newcommand\rank{\mbox{\rm rank}\,}

\newcommand\codim{ \mbox{\rm codim}\, }

\def\qed{\hfill $\Box$}
\def\proof{\noindent {\sl Proof} :\;  }

\def\rd{\partial}

\def\cod{\mbox{{\rm cod}}}
\def\Acod{\A\mbox{{\rm -cod}}\,}

\def\bv{\mbox{\boldmath $v$}}
\def\bw{\mbox{\boldmath $w$}}

\def\bn{\mbox{\boldmath $n$}}

\def\b0{\mbox{\boldmath $0$}}

\def\rk{\mbox{\rm \footnotesize rk}}

\begin{document}

\title[Recognition of plane-to-plane map-germs]
{Recognition of plane-to-plane map-germs}
\author[Y.~Kabata]{Yutaro Kabata}
\address[Y.~Kabata]{Department of Mathematics, 
Graduate School of Science,  Hokkaido University,
Sapporo 060-0810, Japan}
\email{kabata@math.sci.hokudai.ac.jp}
\date{}
%
%

%
\begin{abstract} 
We present a complete set of criteria for determining $\A$-types of 
plane-to-plane map-germs of corank one with $\A$-codimension $\le 6$, which provides   
a new insight into the $\A$-classification theory 
from the viewpoint of recognition problem. 
As an application to generic differential geometry, 
we discuss about projections of smooth surfaces in $3$-space. 
\end{abstract}

\maketitle


\section{Introduction} 
We revisit the $\A$-classification of  local singularities of 
plane-to-plane maps. 
Here $\A$ denotes the group of diffeomorphism germs of source and target planes 
preserving the origin. 
The classification has been achieved by  J. H. Rieger, M. A. S. Ruas \cite{Rieger, Rieger2, RR} 
-- for instance, 
Table \ref{Riegerlist} below shows the list of all corank one map-germs 
with $\A$-codimension $\le 6$. 
When we apply the classification to some specific geometric situation, 
it often becomes a cumbersome task to detect which $\A$-type  a given map-germ belongs to, 
that is referred to as ``$\A$-recognition problem" (cf. \cite{Gaffney}). 
In fact, Rieger's  algorithm frequently uses Mather's Lemma 
to reduce the jet to some nicer form, at which   
the coordinate changes are not explicitly given 
(dotted lines in the recognition trees Fig. 1--5 in \cite{Rieger} 
indicate such processes). To fill up the process is not easy: 
the task is essentially related to  deeper understanding on 
a filtered structure of the $\A$-tangent space of the germ,
as T. Gaffney pointed out in an earlier work \cite{Gaffney}. 

In this paper, we present a complete set of criteria for detecting $\A$-types of corank one germs with 
$\A$-codimension $\le 6$ (Theorem \ref{thm1}). 
That is  a useful package consisting of two-phased criteria, 
which would easily be implemented in computer. 
 The first one is about geometric conditions on `specified jets' for {\it topological $\A$-types
in terms of intrinsic derivatives 
\cite{Whitney, Saji, SUY, Porteous, Levine}, 
and the second is about algebraic conditions on Talyor coefficients of germs with some specified jets, }
which  are obtained by describing explicitly 
all the required coordinate changes of source and target of map-germs 
which are hidden in the classification process 
 (Proposition \ref{nondeg thm}, \ref{singdeg thm}, \ref{doubdeg thm}, and \ref{trideg them}). 


For example, look at  the cases of 
the butterfly $(x,xy+y^5\pm y^7)$ and the elder butterfly $(x,xy+y^5)$,    
which are combined into a single topological $\A$-type. 
Suppose that a map-germ $f=(f_1, f_2): \R^2,0 \to \R^2,0$ with corank one is given. 
Put  $\lambda (x,y):=\frac{\partial (f_1, f_2)}{\partial (x, y)}$, 
and take an arbitrary vector field 
$\eta:=\eta_1(x,y)\frac{\partial}{\partial x}+\eta_2(x,y)\frac{\partial}{\partial y}$ 
near the origin of the source space so that 
$\eta$ spans 
 $\ker df$ on $\lambda=0$. Denote $\eta^k g := \eta(\eta^{k-1} g)$. 
We show that 
the corresponding weighted homogeneous {\it specified jet} (see section \ref{top_A}) 
is characterized in terms of $\lambda$ and $\eta$: 
$$ j^5f(0) \sim_{\A^5} (x,xy+y^5) \;\; \Longleftrightarrow 
\left\{
\begin{array}{l}
d \lambda (0) \neq 0, \\
\eta \lambda (0) = \eta^2 \lambda (0) =\eta^3 \lambda (0) =0, \;\; 
 \eta^4 \lambda (0) \neq 0
\end{array}
\right. 
$$
 Notice that 
the condition in the right hand side does not depend on the choices of local coordinates and 
the null vector field $\eta$.
 The subtle difference between these  ($C^\infty$-)$\A$-types is expressed 
by the following {\it Taylor coefficients condition}: 
If we write $f =(x,xy+y^5+\sum _{i+j\geq 6}a_{ij} x^iy^j)$, 
then 
$$\textstyle 
f \sim_\A (x,xy+y^5\pm y^7) \;\; \Longleftrightarrow \;\; 
a_{07}-\frac58 a_{06}^2 \not=0,$$
otherwise, $f$ is of type elder butterfly. 
It should be noted that  for the butterfly
T. Gaffney \cite{Gaffney} found the same condition  on Taylor coefficients 
by studying the structure of $\A$-tangent space (Example 1.4. in \cite{Gaffney}). 
Our approach is more direct by extending the method used in \cite{Rieger, Bruce}, 
 and we describe such conditions for all $\A$-types in Rieger's list ($\A$-codimension $\le 6$).  

 Our second purpose is to demonstrate 
a systematic use of our criteria 
for map-germs arising in some specific geometric situation. 
We develop a method of J. W. Bruce \cite{Bruce} 
for an application to {\it extrinsic differential geometry} of surfaces. 
Look at a generic surface in $\R^3$ from a viewpoint (camera), 
then we get locally a smooth map from the surface to the plane (screen), that is called 
the {\it central projection}. Their singularities have been classified by V. I. Arnold and O. A. Platonova 
(also O. P. Shcherbak, V. V. Goryunov) \cite{Arnold, Arnold2, Goryunov, Platonova, Shcherbak} 
based on a different framework. 
It is shown that 
some germs of $\A$-codimension $5$ {\it do not appear} generically in central projections, 
although the reason has not been quite clear from the context of $\A$-classification,  
as Rieger noted in his paper \cite{Rieger}. 
Our criteria make the reason very clear --   
the condition of intrinsic derivatives $\eta^k\lambda$ determines  jets of Monge form of the surface, 
while the condition of  Taylor coefficients determines a special position of viewpoints 
(Remarks \ref{h-focal} and \ref{p-focal}). 
We present an alternative transparent proof of Arnold-Platonova's theorem 
within the $\A$-classification theory, moreover, 
we classify singularities arising in central projections of 
{\it moving surfaces} with one-parameter in $3$-space (Theorem \ref{mainthm2}).




 As a byproduct, in another paper \cite{SKDO} 
we obtain a generalization of projective classification  
of jets of Monge forms by Platonova \cite{Platonova}. 
Our criteria are also useful to determine the bifurcation diagrams of map-germs, 
especially of corank two. See \cite{YKO, YKO2} for the detail.

The rest of this paper is organized as follows. 
In \S 2 we briefly introduce the classification of plane-to-plane map-germs. 
In \S 3 we give a complete set of criteria for all $\A$-types with $\Acod \le 6$.
 In \S 4 we show an application of our criteria to the central projection of smooth surfaces.

\section{Preliminary}
\subsection{$\A$-classification}
To begin with, we briefly summarize the basics of singularity theory of map-germs. 
Let $\E_n$ be an $\R$-algebra of smooth map-germs $\R^n, 0\rightarrow \R$ 
with a unique maximal ideal $m_n$. 
The $\E_n$-module consisting of map-germs $\R^n, 0\rightarrow \R^p, 0$ is isomorphic to 
$m_n \E_n^p$. On this space, 
the group of diffeomorphism germs defines an equivalence relation: 
$f, g : \R^n, 0\rightarrow \R^p, 0$ are {\it $\A$-equivalent} ($f\sim_{\A}g$), 
if there exist diffeomorphism germs $\phi$ and $\psi$ of $\R^n,0$ and $\R^m,0$ so that $f=\psi\circ g \circ \phi^{-1}$.
We denote by $\A.f$ the $\A$-orbit of $f$. 
If an $\A$-orbit has finitely many nearby orbits, then the orbit is called {\it $\A$-simple}; 
otherwise, there is some family of $\A$-orbits,  called an {\it $\A$-moduli}. 

Let $\xi : \R^n, 0 \rightarrow T\R^p$ be a smooth map-germ such that 
$\pi \circ \xi =f$ (where $\pi$ is a projection of tangent vector bundle). 
We call $\xi$ the {\it the vector field along $f$} 
or {\it infinitesimal deformation of $f$}, 
and denote the set of all the the vector field along $f$ by $\theta(f)$. 
In an obvious way, $\theta(f)$ is a $\E_n$-module. 
For the identity maps $id_n : \R^n, 0 \rightarrow \R^n, 0$,  $id_p : \R^p, 0 \rightarrow \R^p, 0$, 
we write $\theta(n)=\theta(id_n)$, $\theta(p)=\theta(id_p)$, 
which are the module of vector field-germs. 
We define $tf : \theta(n) \rightarrow \theta(f)$ by the map $\xi \mapsto df \circ \xi$, and $\omega f : \theta(p) \rightarrow \theta(f)$ by the map $\eta \mapsto \eta \circ f$.
With these notations above, we define the {\it $\A$-tangent space of $f$} by 
$$T\A(f):=tf(m_n \theta_n)+\omega f(m_p \theta_p) \; \subset m_n \theta(f).$$
In fact, this space consists of all vectors 
$\frac{d}{dt} (\psi_t \circ f \circ \phi_t^{-1}) \big|_{t=0}$ 
where $\psi_t$ and $\phi_t$ are deformations of identity maps with 
$\psi_t(0)=0$, $\phi_t(0)=0$. 
We define the {\it $\A$-codimension}  ($\Acod$) of $f$ by 
$\cod(\A, f):=\mbox{dim}_{\R} \, m_n \theta(f)/T\A(f)$. 

We define {\it the $r$-jet space} $J^r(n, p)$ to be the set of $r$-jets of map-germs at the origin. 
This is naturally identified with $m_n\E^p_n/m^{r+1}_n\E^p_n$. 
We say $j^rf(0)$ {\it is $\A^r$-equivalent to} $j^rg(0)$ 
($j^rf(0)\sim_{\A^r}j^rg(0)$) 
if there exist diffeomorphism germs 
$\psi:\R^n, 0 \rightarrow \R^n, 0$ and $\phi:\R^p, 0\rightarrow \R^p, 0$ 
such that 
$j^rf(0)=j^r(\phi\circ g \circ \psi^{-1}) (0)$. 
The $r$-jets of diffeomorphism germs form a Lie group denoted by 
$\A^r$, which naturally has an algebraic action on $J^r(n, p)$, 
hence 
 the orbit $\A^r(j^rf(0))$ is a locally closed semi-algebraic submanifold.  
For a smooth map-germ $f : \R^n, 0\rightarrow \R^p, 0 $,
we say $f$ is {\it $r$-$\A$-determined}, 
 if $f\sim_{\A}g$ holds for any smooth map-germs $g : \R^n, 0\rightarrow \R^p, 0 $ such that
 $j^rf(0)=j^rg(0)$. 
 When $f$ is $r$-determined for some $r$, 
 we say $f$  {\it is finitely $\A$-determined}.

From now on, we consider the case of $n=p=2$. 
We are concerned with 
the $\A$-classification  \cite{Rieger, Rieger2, Rieger3, RR}: In particular, 
all $\A$-types of $f: \R^2,0 \rightarrow \R^2,0$ of corank one with $\A$-cod at most 6 
are listed in Table \ref{Riegerlist}.   
Here we use the notation $\Acod$ for germs with moduli to refer to the codimension of stratum. 
There are 29 types in the list with additional sign $\pm$, and 
we use Rieger's notation $1, 2, \cdots, 19$ for the $\A$-types throughout this paper. 
The type no.14: $(x, xy^2+y^5)$ is not included in Table \ref{Riegerlist}, 
since it has $\A$-codimension $7$.

\begin{table}\label{Riegertable1}
$$
\begin{array}{c| l | l }
\mbox{$\A$ -cod}  &\mbox{type}  & \mbox{normal form}\\
\hline
0&1\, (\mbox{regular}) &  (x,y) \\
\hline
1&2\, (\mbox{fold}) &  (x,y^2) \\
\hline
2&3\, (\mbox{cusp}) & (x,xy+y^3) \\
\hline
3&4_2\, (\mbox{beaks and lips}) &(x, y^3\pm x^2y) \\
&5(\mbox{swallowtail}) & (x,xy+y^4) \\
\hline
4&4_3\, (\mbox{goose}) &(x, y^3+ x^3y) \\
&6\, (\mbox{butterfly}) & (x,xy+y^5\pm y^7) \\
&11_5\, (\mbox{gulls}) & (x,xy^2+y^4+y^5) \\
\hline
5&4_4\, (\mbox{ugly goose}) & (x,y^3 \pm x^4y) \\
&7\, (\mbox{elder butterfly}) & (x,xy+y^5) \\
&11_7\, (\mbox{ugly gulls}) & (x,xy^2+y^4+y^7) \\
&12 & (x,xy^2+y^5+y^6) \\
&16 & (x,x^2y+y^4\pm y^5) \\
& 8\, \mbox{(unimodal)} &  (x,xy+y^6 \pm y^8+\alpha y^9) \\
\hline
6&4_5 & (x,y^3 + x^5y) \\
&9& (x,xy+y^6+y^9) \\
&10^\dagger\,  \mbox{(bimodal)} & (x,xy+y^7\pm y^9 +\alpha y^{10} + \beta y^{11}) \\
&11_9& (x,xy^2+y^4+y^9) \\
&13& (x,xy^2+ y^5\pm y^9) \\
&15\,  \mbox{(unimodal)} & (x,xy^2+y^6+y^7+\alpha y^9) \\
&17& (x,x^2y+y^4) \\
&18^\dagger\,  \mbox{(bimodal)} & (x,x^2y+xy^3+\alpha y^5+y^6 + \beta y^7) \\
&19\,  \mbox{(unimodal)} & (x,x^3y+\alpha x^2y^2+y^4 + x^3y^2) \\
\end{array}
$$
\caption{$\A$-classification up to $\A$-$\mbox{cod}\le 6$ \cite{Rieger}. 
$\dagger$: excluding exceptional values of the moduli}
\label{Riegerlist}
\end{table}

\subsection{Topological $\A$-classification} \label{top_A}
Two germs are {\it topologically $\A$-equivalent} if they commute 
via some homeomorphisms of source and target;   that is 
the version where one just replaces diffeomorphisms for $\A$-equivalence by homeomorphisms. 
By using a theorem of J. Damon \cite{Damon}, 
several different $\A$-types in Table \ref{Riegerlist} 
are combined into a single topological $\A$-type: 
Those are listed in the following Table \ref{Riegerlist_top} \cite{Rieger2}. 

\begin{table}
$$
\begin{array}{c | c | l }
\mbox{topoligical type} & \mbox{$\A$ -type}   & \mbox{normal form}\\
\hline 
I_k^\pm \; (k\ge 2) & 4_k^\pm & (x, y^3\pm x^ky) \\
II_k \; (4\le k \le 6) & 5-10 & (x, xy+y^k)\\
III_k\;  (k \ge 2) & 11_{2k+1} & (x, xy^2+y^4+y^{2k+1})\\
IV_5 & 12, 13, (14) & (x, xy^2+y^5) \\
V_1 & 16, 17 & (x, x^2y+y^4)
\end{array}
$$
\caption{Some different $\A$-types (with $\A\mbox{-codim}\, \le 6$) are combined into 
the same topological $\A$-types \cite{Rieger2}.}
\label{Riegerlist_top}
\end{table}

We introduce a coarser classification than topological $\A$-classification 
for our convenience. 
We provisionally call 
the weighted homogeneous part of each normal form 
in Table \ref{Riegerlist_top} 
the {\it specified jet} for the corresponding topological $\A$-type, 
except for $4_k$-types; 
the specified jet of $4_k$ ($k \ge 3$) 
is defined to be $(x, y^3)$.  
Note that both germs $(x, y^3)$  and $(x, xy^2+y^4)$ 
are not finitely $\A$-determined, 
thus we can not use Damon's theorem \cite{Damon};   
indeed $4_k$ and $11_{2k+1}$ for different $k$ may have different topological $\A$-types. 
However it is useful for our purpose 
to gather all $4_k$  of  $k \ge 3$ (resp. $11_k$) into 
a group  $I_*$ (resp. $III_*$) of $\A$-types having the same specified jet. 

Here we list up all specified jets of germs under consideration 
in this paper (stable germs are omitted and 
specified jets  of types $15, 18, 19$ are denoted by $IV_6$, $V_2$, $VI$, respectively): 
\begin{eqnarray*}
&& I_2:  (x, y^3\pm x^2y), \quad I_*: (x, y^3),\\
&& II_4: (x, xy+y^4), \quad  II_5: (x, xy+y^5), \quad  II_6: (x, xy+y^6),\quad II_7: (x, xy+y^7)\\
&& III_*: (x, xy^2+y^4)\\
&& IV_5: (x, xy^2+y^5), \quad IV_6: (x, xy^2+y^6),\\
&& V_1: (x, x^2y+y^4), \quad V_2: (x, x^2y+xy^3) \;  \mbox{or} \; (x, x^2y)\\
&& VI: (x, y^4+\alpha x^2y^2+x^3y)\;  \mbox{or} \; (x, y^4+\alpha x^2y^2)
\end{eqnarray*}


\section{Criteria for map-germs}

We state our main result: 

\begin{thm}\label{thm1}
Specified jets 
of topologically $\A$-equivalent types of plane-to-plane germs 
 with $\A$-codimension up to $6$ are explicitly characterized   
by means of geometric terms $\lambda$ and $\eta$ 
as in Table \ref{criteria1}:   Precisely saying, 
given a map-germ $f$  of corank one, the jet $j^rf(0)$ is $\A^r$-equivalent to one of 
the specified $r$-jets listed in Table \ref{criteria1} if and only if 
the corresponding condition of $\lambda$ and $\eta$ for $f$ in Table \ref{criteria1}
 is satisfied. 
A complete set of criteria for detecting $\A$-types of germs with $\A$-codimension up to $6$ 
(Table \ref{Riegerlist}) is achieved by  
 adding conditions in coefficients of Taylor expansions,  
which are precisely described    
in Proposition \ref{nondeg thm}, \ref{singdeg thm}, \ref{doubdeg thm}, and 
\ref{trideg them} below. 
\end{thm}

\begin{table}
{\small 
$$
\begin{array}{ l | l | l | l }
 \mbox{specified jet }  &\mbox{$\A$-type}  & \mbox{condition}&\mbox{type}\\
\hline
\hline
II_4:  (x, xy+y^4)&  5 &
d\lambda(0)\not=0, & \\
&& \eta\lambda(0)=\eta^2\lambda(0)=0,&\;\; A_0\;/\;4 \\
&& \eta^3\lambda(0)\not=0&\\
 \hline
II_5:  (x, xy+y^5)&  6, 7 &
d\lambda(0)\not=0, &\\
&& \eta\lambda(0)=\eta^2\lambda(0)=\eta^3\lambda(0)=0, &\;\; A_0\;/\;5\\
&& \eta^4\lambda(0)\not=0 &\\
\hline
 II_6: (x, xy+y^6)&  8, 9 &
 d\lambda(0)\not=0, &\\
&&
\eta\lambda(0)=\cdots=\eta^4\lambda(0)=0&\;\; A_0\;/\;6\\
&& \eta^5\lambda(0)\not=0&\\
\hline
II_7: (x, xy+y^7)&  10 &
 d\lambda(0)\not=0, & \\
&&
\eta\lambda(0)=\cdots=\eta^5\lambda(0)=0,& \;\; A_0\;/\;7  \\
&&\eta^6\lambda(0)\not=0&\\
\hline
\hline
I_2:  (x, y^3\pm x^2y)& 4_2^\pm  &
d \lambda (0) = 0, &\\
&&\det H_\lambda (0) \not=0,&\;\; A_1\; / \;3\\ 
&& {\eta} ^2 \lambda (0) \neq 0&\\
\hline
III_*: (x, xy^2+y^4)& 11_{odd} &
 d\lambda(0)=0, \det H_\lambda(0)<0, &\\
&&\eta^2\lambda(0)=0,  &\;\; A^{-}_1\; / \;4\\
&& \eta^3\lambda(0)\not=0 &\\
\hline
IV_5: (x, xy^2+y^5)& 12, 13 &
d\lambda(0)=0, \det H_\lambda(0)<0,&\\
&&  \eta^2\lambda(0)=\eta^3\lambda(0)=0, &\;\; A^{-}_1\; / \;5\\
&& \eta^4\lambda(0)\not=0&\\
\hline
IV_6: (x, xy^2+y^6)& 15 &
d\lambda(0)=0, \det H_\lambda(0)<0,&\\
&&  \eta^2\lambda(0)=\eta^3\lambda(0)=\eta^4 \lambda(0)=0, & \;\; A^{-}_1\; / \;6\\
&&\eta^5\lambda(0)\not=0&\\
\hline
\hline
I_*:  (x, y^3)& 4_* &
d \lambda (0) = 0, &\\
&&\rk H_\lambda (0) =1, &\;\; A_k\; /  \; 3 \\ 
&& {\eta} ^2 \lambda (0) \neq 0 &\\
\hline
V_1: (x, x^2y+y^4)& 16, 17 &
d\lambda(0)=0, \rk H_\lambda(0)=1, &\\
&&\eta^2\lambda(0)=0, & \;\; A_2\; /  \; 4\\
&&\eta^3\lambda(0)\not=0 &\\
\hline
V_2: (x, x^2y+xy^3),& 18 &
d \lambda (0) = 0, &\\
\;\;\; (x, x^2y)&&\rk H_\lambda (0) =1, &\;\; A_k \; / \ge 5 \\ 
&&
{\eta} ^2 \lambda (0) ={\eta} ^3 \lambda (0) =0 &\\
\hline
\hline
VI: (x,y^4+\alpha x^2y^2),& 19 &
d \lambda (0) = 0, &\\
\; (x,y^4+\alpha x^2y^2+x^3y)&& \rk  H_\lambda (0) =0,& \;\; D_4\; / \;4  \\
&&\eta^3\lambda(0) \neq 0 &\\
\hline
\end{array}
$$
}
\caption{Criteria for  plane-to-plane germs with $\A$-codimension up to $6$ 
  (stable germs are omitted).  Here 
 the last column means singularity types of $\lambda$ at $0$ 
 and local degree of complexified germs respectively. 
 Refer to Rieger's original list (Table 1 in \cite{Rieger}) for other geometrical invariants. }
\label{criteria1}
\end{table}

\begin{rem}\label{coeff}
{ \rm 
Our condition in coefficients of Taylor expansions detects $\A$-types
among types having the same specified jets, 
however the geometric meaning is not so clear. 
For a few cases, 
T. Gaffney \cite{Gaffney} found the same conditions in studying 
a finer algebraic structure of the corresponding $\A$-tangent space. 
It would be interesting to compare these two approaches. 
It would also be reasonable to discuss about the problem in the context of Damon's
$\mathcal{K}_D$-theory using the logarithmic vector fields along 
the $A_\mu$-type discriminant of a stable unfolding. 
That will be considered somewhere else.}
\end{rem}


The proof is divided into the following four cases: 

\

\begin{tabular}{ll}
({\it case 0}) & $d\lambda(0) \not=0$; \\
({\it case 1}) & $d\lambda(0)=0$ and  $\rk H_\lambda (0)=2$; \\ 
({\it case 2}) & $d\lambda(0)=0$ and  $\rk H_\lambda (0)=1$;\\
({\it case 3}) & $d\lambda(0)=0$ and  $\rk H_\lambda (0)=0$.  
\end{tabular}

\

In fact, 
Table \ref{criteria1} is separated into these four cases by double lines. 
These cases deal with the same process in recognition trees Fig. 1--5 
in \cite{Rieger}: 
Cases 0, 1, 3 correspond to Fig. 1, 3, 5, respectively, 
and case 2 corresponds to both Fig. 2 and 4 in \cite{Rieger}.

For our simplicity, we omit the 
case of $\Acod \le 3$, that is the set of characterizations by Whitney and Saji in \cite{Saji, SUY, Whitney}. 
In the following proof, we frequently use 
Rieger's results (e.g., $\A$-determinacy of germs), 
which should be referred to \cite{Rieger}.

\subsection{{\bf Case 0}: $d\lambda \not= 0$  ($S(f)$ is smooth)}

 We deal with types $6-10$ of $\Acod =  4, 5, 6$. 

\begin{prop}\label{nondeg thm}
 For a plane-to-plane map-germ $f$ of corank one, 
\begin{enumerate}
\renewcommand{\labelenumi}{(\arabic{enumi})}
\item
For $r\ge5$,\\
$j^rf(0) \sim_{\A^r} (x,xy+y^r)$ $\Longleftrightarrow$ 
$$
d \lambda (0) \neq 0, \;\; \eta^{i} \lambda (0) =0\; (1\le i \le r-2), \;\; \eta^{r-1} \lambda(0)\ne0.
$$

\item
 If we write 
$f =(x,xy+y^5+\sum _{i+j\geq 6}a_{ij} x^iy^j)$, 
\begin{eqnarray*}
\textstyle a_{07}-\frac58 a_{06}^2\not=0  &\Longleftrightarrow &  f \sim_{\A} (x,xy+y^5\pm y^7)\cdots \fbox{6}\,, \\
\textstyle a_{07}-\frac58 a_{06}^2=0  &\Longleftrightarrow &  f \sim_{\A}(x,xy+y^5)\cdots \fbox{7}\,. 
\end{eqnarray*}

\item 
 If we write 
$f =(x,xy+y^6+\sum _{i+j\geq 7}a_{ij} x^iy^j)$, 
$$\textstyle
\begin{array}{lcc}
\vspace{2mm}
 a_{08}-\frac{3}{5} a_{07}^2 \not=0  &\Longleftrightarrow&   f \sim_{\A} (x,xy+y^6\pm y^8+\alpha y^9)\cdots \fbox{8}\,, \\
\left\{\begin{array}{l}
\vspace{1mm}
a_{08}-\frac{3}{5} a_{07}^2=0  \\
{a}_{09}-\frac{7}{25} a_{07}^3\neq 0 
\end{array}\right.  &\Longleftrightarrow&  f \sim_{\A}(x,xy+y^6+y^{9})\cdots \fbox{\mbox{$9$}}\,.
\end{array}
$$

\item 
If we write 
$f =(x,xy+y^7+\sum _{i+j\geq 8}a_{ij} x^iy^j)$, 
\begin{eqnarray*}
\textstyle a_{09}-\frac{7}{12} a_{08}^2\not=0 &\Longleftrightarrow &  j^{11}f(0) \sim_{\A^{11}} (x,xy+y^7\pm y^9+\alpha y^{10}+\beta y^{11}),
\end{eqnarray*}
 and excluding exceptional values of $\alpha$ and $\beta$
 $$
 f \sim_{\A} (x,xy+y^7\pm y^9+\alpha y^{10}+\beta y^{11})\cdots \fbox{10}\,. 
 $$

\end{enumerate}
\end{prop}

In order to prove 1 in Proposition \ref{nondeg thm}, we need the next lemma 
based on Lemma $2.6$ in \cite{Saji}. 
Notice that $\lambda$ is changed by multiplying a non-zero function when we take another coordinates, 
and also that there is an ambiguity to choose the null vector field $\eta$. 

\begin{lem}\label{nondeg invariant lem}
The conditions on the right hand side of 
 1 in Proposition \ref{nondeg thm} 
are independent from the choice of coordinates of the source and target 
and the choice of $\eta$. 
\end{lem}
{ \proof The proof is similar to that of Lemma $2.6$ in \cite{Saji}. \qed}

\

\noindent {\sl Proof of 1 in Proposition \ref{nondeg thm}} :\;  
It is easily checked that 
for the $r$-jet $(x,xy+y^r)$, the condition in the right hand side holds.
Thus the ``only if'' part of 1 follows from Lemma \ref{nondeg invariant lem}.  

The ``If" part is shown by finding a suitable coordinate change. 
Assume that the condition on the right hand side of 1 holds for $f$. 
Since $f$ is of corank $1$ at $0$, 
we may write 
$$\textstyle f(x,y)=(x, \sum_{ i+j \ge 2} a_{ij} x^i y^j).$$
Take $\eta=\frac{\partial}{\partial y}$ 
and 
$$\textstyle 
\lambda(x,y)=\sum_{i+j \ge 2} j\cdot a_{ij} x^i y^{j-1}.
$$
For this choice of coordinates and $\eta$, 
by Lemma \ref{nondeg invariant lem}, we have 
$$
 d \lambda (0) \neq 0, \;\; \eta \lambda (0) = \eta^2 \lambda (0) =\cdots=\eta^{r-2} \lambda (0) =0, \;\; 
 \eta^{r-1} \lambda (0) \neq 0. 
$$
Then 
$a_{11}\neq 0$, $a_{02}=a_{03}=\cdots=a_{0\, r-1}=0$, $a_{0r}\neq 0$.
By some coordinate change, 
we have  $j^rf(0)=
(x,xy+y^r)$. 
\qed

\

 The following proof of the claim 2 uses a simple trick for eliminating a certain term 
 in the normal form.  { This trick is standard for the classification of map-germs as seen in Bruce's work \cite{Bruce}, and will implicitly appear several times in other cases. }

\vspace{1mm}

\noindent {\sl Proof of 2 in Proposition \ref{nondeg thm}} :\; 
{ Since both 6-type and 7-type are $7$-determined, 
our task is to show} 
$$\textstyle
(x, xy+y^5+cy^6+dy^7) \sim_{\A^7} (x, xy+y^5+(d-\frac{5}{8} c^2)y^7). 
$$

Write $xy+y^5+cy^6+dy^7=xy+y^5(1+\alpha y)+\beta y^6+dy^7$ 
with $\alpha+\beta=c$. 
By the coordinate change so that 
$$\bar{x}=x, \quad \bar{y}^5=y^5(1+\alpha y), $$
the $7$-jet has the form
$$
\textstyle(\bar{x}, \bar{x}(\bar{y}-\frac15 \alpha \bar{y}^2+\frac{4}{25} \alpha^2 \bar{y}^3+h.o.t.)+\bar{y}^5+\beta (\bar{y}-\frac15 \alpha \bar{y}^2+h.o.t.)^6 +d \bar{y}^7).  
$$
By the coordinate change  
\begin{eqnarray*}
&&\textstyle x=
\bar{x}(1-\frac15 \alpha \bar{y}+\frac{4}{25} \alpha^2 \bar{y}^2+h.o.t.)\\
&&y=\bar{y}
\end{eqnarray*}
 the jet is written by 
$$\textstyle
(x(1+\frac15 \alpha y-\frac{3}{25} \alpha^2 y^2+h.o.t.) , \; 
xy+y^5+\beta y^6 +(d-\frac65 \alpha \beta) y^7). 
$$
Then, by the coordinate change of target
$$\textstyle 
(X,Y)\rightarrow (X-\frac15 \alpha Y, Y),$$ 
we eliminate the term $\frac{1}{5}\alpha xy$ in the first component; 
Hence the jet becomes 
$$
\left(
\begin{array}{l}  
 x(1-\frac{3}{25} \alpha^2 y^2+h.o.t.)-\frac15 \alpha(y^5+\beta y^6+(d-\frac65 \alpha \beta)y^7) \\
xy+y^5+\beta y^6 +(d-\frac65 \alpha \beta) y^7
 \end{array}
 \right).
$$
Take $\tilde{x}$ to be the first component 
and $\tilde{y}=y$, 
then  the jet is written by (after rewriting variables)
$$\textstyle
 (x , xy(1+\frac{3}{25}\alpha^2y^2+h.o.t.)+y^5+(\frac{1}{5}\alpha + \beta)y^6+(d- \alpha \beta) y^7).  
$$
Now we choose $\alpha=\frac54 c$ and $\beta=-\frac14c$ 
to kill the term $y^6$ in the second component. 
Finally by 
$\tilde{x}= x$ and 
$\tilde{y}=y(1+\frac{3}{25}\alpha^2y^2+h.o.t.)$,   
we obtain the form 
$$\textstyle
 (x, xy+y^5+(d-\frac{5}{8} c^2)y^7).
$$
This completes the proof. 
\qed

\

\noindent {\sl Proof of 3 in Proposition \ref{nondeg thm}} :\; 
The proof is similar to that of the claim 2 just described above: 
First we eliminate the terms including $x$ of order $\ge 7$, and then 
 we directly show that 
$$ \textstyle
(x, xy+y^6+cy^7+dy^8+ey^9) \sim_{\A^9} (x, xy+y^6+(d-\frac35 c^2)y^8+(e-\frac{7}{5}cd+\frac{14}{25}c^3)y^9).
$$
In fact, rewriting variables as $\tilde{x}, \tilde{y}$ of the germ in the left hand side,  
substitute 
\begin{eqnarray*}
&&\textstyle
\tilde{x}=x + \frac {c } {5} x y + \frac {c } {5}y^6 - \frac {3 c^3 } {25} y^8 +\frac {c d} {5}  y^8+\frac{14c^4}{125}y^9-\frac{7c^2d}{25}y^9+\frac{ce}{5}y^9, \\
&&\textstyle
\tilde{y}=y - \frac {c 
} {5} y^2 + \frac {c^2 } {25} y^3 - \frac {c^3 } {125} y^4+ \frac {c^4 
} {625} y^5 - \frac {c^5 } {3125} y^6\\
&&\;\;\;\;\;\;\;\;\;\;\;\;\textstyle + \frac {c^6 } {15625}y^7-\frac{c^7}{78125}y^8+\frac{2c^8}{390625}y^9, 
\end{eqnarray*}
and take the coordinate change of the target 
$$\textstyle(X,Y)\mapsto \left(X-\frac{c}{5} Y,Y \right), $$
then we get the equivalence. 
Here $c=a_{07}$, $d=a_{08}$, $e=a_{09}$, and 
both $8$-type and $9$-type are $9$-determined, thus we have the claim 3. \qed

\

\noindent {\sl Proof of 4 in Proposition \ref{nondeg thm}} :\; 
Also in a similar way as above we see 
$$ \textstyle 
(x, xy+y^7+cy^8+dy^9) \sim_{\A^9} (x, xy+y^7+(d-\frac{7}{12} c^2)y^9). 
$$
In fact, it is achieved by 
\begin{eqnarray*}
&&\textstyle 
\tilde{x}=
x + \frac {c } {6}x y + \frac {c } {6}y^7 - \frac {7 c^3 } {72}y^9 + 
\frac {c d} {6}  y^9 , \\
&&\textstyle
\tilde{y}=y - \frac {c } {6}y^2 + \frac {c^2} {36} y^3 - \frac {c^3 } {216}y^4 
+ \frac {c^4 } {1296}y^5 - \frac {c^5 } {7776}y^6 + \frac {c^6 } 
{46656}y^7 \\
&&\textstyle \qquad \quad  
- \frac {c^7 } {279936}y^8 - \frac {5 c^8 } {93312}y^9
\end{eqnarray*}
and 
$$\textstyle 
(X,Y)\mapsto (X-\frac{c}{6}Y, Y).$$ 
In addtion, we easily get 
$$\textstyle
(x, xy+y^7\pm y^9+O(10))\sim_{\A} (x,xy+y^7\pm y^9+\alpha y^{10}+\beta y^{11}+O(12))
$$
for some $\alpha,\,\beta\in\R$. 
In Rieger \cite{Rieger}, it is shown that 
the $10$-type is $11$-determined 
for generic $\alpha$ and $\beta$ excluding some values 
explicitly given in \cite[p.359]{Rieger}. 
This implies the claim 4. 
\qed

\subsection{{\bf Case 1}:  $d\lambda(0)=0$,   $\rk H_\lambda (0)=2$} 

 We deal with types $11_{2k+1} \, (k=2,3,4)$, \,$12$,\,$13$,\,$15$ of $\Acod =4, 5, 6$.

\begin{prop}\label{singdeg thm}
 For a plane-to-plane map-germ $f$ of corank one, 
\begin{enumerate}
\renewcommand{\labelenumi}{(\arabic{enumi})}
\item
For $r\ge4$,\\
$j^r f(0) \sim_{\A^r} (x,xy^2+y^r)$ $\Longleftrightarrow$ 
$$
d \lambda (0) = 0, \;\; 
\det H_\lambda (0)< 0, \;\; \eta^{i} \lambda (0) =0\; (2\le i \le r-2), \;\; \eta^{r-1} \lambda(0)\ne0.
$$

\item
If we write 
$f = (x,xy^2+y^4+\sum _{i+j\geq 5}a_{ij} x^iy^j)$, 
$$\textstyle
\begin{array}{ccl}
 a_{05}\neq 0  &\Longleftrightarrow&   f \sim_{\A} (x,xy^2+y^4+y^{5})\cdots \fbox{\mbox{$11_5$}}\,, \\
\left\{\begin{array}{l}
a_{05}=0  \\
{a}_{07}-2a_{15}+4 a_{23}\neq 0
\end{array}\right.  &\Longleftrightarrow&  f \sim_{\A}(x,xy^2+y^4+y^{7})\cdots \fbox{\mbox{$11_7$}}\,, \\
a_{05}={a}_{07}-2a_{15}+4 a_{23}=0 &\Longleftrightarrow& j^7f(0) \sim_{\A^7}(x,xy^2+y^4).
\end{array}
$$
Furthermore, if we write 
$f = (x,xy^2+y^4+\sum _{i+j\geq 8}c_{ij} x^iy^j)$, 
$$\textstyle
c_{09}-2c_{17}\neq 0
 \Longleftrightarrow  f \sim_{\A}(x,xy^2+y^4+y^{9})\cdots \fbox{\mbox{$11_9$}}\,.
$$

\item
If we write 
$
  f =(x,xy^2+y^5+\sum _{i+j\geq 6}a_{ij} x^iy^j) ,
$ 
$$\textstyle
\begin{array}{ccl}
 a_{06}\neq 0  &\Longleftrightarrow&  f \sim_{\A}(x,xy^2+y^5+y^{6})\cdots \fbox{\mbox{$12$}}\,, \\
\left\{\begin{array}{l}
a_{06}=0  \\
a_{09}-\frac52 a_{16}-\frac56 a_{07}^2\not=0
\end{array}\right.  &\Longleftrightarrow&  f \sim_{\A}(x,xy^2+y^5+y^9)\cdots \fbox{$\mbox{13}$}\,. \\
\end{array}
$$

\item
If we write 
$
  f =(x,xy^2+y^6+\sum _{i+j\geq 7}a_{ij} x^iy^j) ,
$ 
\begin{eqnarray*}
\textstyle a_{07}\neq 0  &\Longleftrightarrow&  f \sim_{\A}(x,xy^2+y^6+y^7+\alpha y^9)\cdots \fbox{15}\,. \\
\end{eqnarray*}

\end{enumerate}
\end{prop}

Note that in claim 1 of Proposition \ref{singdeg thm}, $d \lambda(0)=0$ implies $\eta \lambda(0)=0$. 


\

\noindent{\it Proof of 1 in Proposition \ref{singdeg thm}} :\;
The proof is similar to that of 1 in Proposition \ref{nondeg thm}. 
\qed

\

\noindent {\sl Proof of 2 in Proposition \ref{singdeg thm}} :\; 
Let 
$ f=(x,xy^2+y^4+\sum _{i+j\geq 5}a_{ij} x^iy^j).$
{ By routine coordinate changes, $f$ is equivalent to}
$$\textstyle 
(x,xy^2+y^4+{a}_{05}y^5+ \sum _{i+j\geq 6}{a}'_{ij} x^iy^j)
$$
for some ${a}'_{ij}$ (${a}'_{06}=a_{06}-2a_{14}$ etc). 
Since $(x, xy^2+y^4+y^5)$ is 5-determined,
$a_{05} \not=0$ leads to
$$\textstyle 
f \sim_{\A}
(x, xy^2+y^4+y^5). 
$$

Next we suppose $a_{05}=0$.
A similar coordinate change as above shows 
that $f$ is equivalent to 
$$\textstyle 
(x,xy^2+y^4+(a_{06}-2a_{14})y^6+({a}_{07}-2a_{15}+4 a_{23})y^7+ \sum _{i+j\geq 8}{b}_{ij} x^iy^j) 
$$
for some $b_{ij}$.
Since the germ  $(x,xy^2+y^4+y^7)$ is $7$-determined,
we want to eliminate the term $y^6$ from the second component of the right-hand side.
By a similar argument as in the proof of 2 in Proposition \ref{nondeg thm}, 
we have  
$$
(x, xy^2+y^4+cy^6+dy^7)\sim_{\A^7}(x, xy^2+y^4+dy^7). 
$$
In fact this is achieved by an explicit coordinate change of source 
(writing variables as $\tilde{x}, \tilde{y}$ of the germ in the left hand side) 
\begin{eqnarray*}
&&\textstyle
\tilde{x}={x} + c {x} {y}^2 + c {y}^4 + 
 c d {y}^7\\
&&\textstyle 
\tilde{y}={y} - \frac {c}{2} {y}^3 
+ \frac {3c^2}{8}{y}^5 -
 \frac {9c^3}{16}{y}^7
\end{eqnarray*}
and $(X,Y)\mapsto (X-cY,Y)$ of target. 
With this coordinate change and some adding coordinate change of source, $f$ is equivalent to
$$\textstyle 
(x,xy^2+y^4+({a}_{07}-2a_{15}+4 a_{23})y^7+\sum _{i+j\geq 8}b'_{ij}x^iy^j). 
$$
for some $b'_{ij}$.
Hence ${a}_{07}-2a_{15}+4 a_{23}\not=0$ leads to 
$$\textstyle 
f \sim_{\A}
(x, xy^2+y^4+y^7). 
$$

Finally suppose ${a}_{07}-2a_{15}+4 a_{23}=0$; Put  
$f=(x,xy^2+y^4+\sum _{i+j\geq 8}c_{ij}x^iy^j)$
for some $c_{ij}$. We see 
$$ 
(x, xy^2+y^4+cy^8+dy^9) \sim_{\A^9} (x, xy^2+y^4+dy^9) 
$$
by the change of  source 
(writing variables as $\tilde{x}, \tilde{y}$ of the germ in the left hand side)  
\begin{eqnarray*}
&&\textstyle
\tilde{x}=
x - \frac {c} {2}  x^2 y^2 - \frac {c} {2}  x y^4 + \frac {c^2} {4} 
 x^3 y^4 + \frac {c^2} {2}  x^2 y^6 + \frac {c^2} {4}  x y^8\\
&&\textstyle
\tilde{y}=y + \frac {c} {4}  x y^3 - \frac {c}{4} y^5- \frac {c^2} {32}  
x^2 y^5 - \frac {c^2} {16}  x y^7 - \frac {5 c^2} {8} y^9
\end{eqnarray*}
and the change of target 
$$\textstyle (X,Y)\mapsto (X+\frac{c}{2} X Y, Y).$$
Then it turns out that 
$f$ is $\A$-equivalent to 
$$\textstyle 
(x,xy^2+y^4+(c_{09}-2c_{17})y^9+O(10)).
$$
Since $(x, xy^2+y^4+y^9)$ is $9$-determined, $c_{09}-2c_{17}\not=0$ 
leads to 
$$\textstyle 
f \sim_{\A}
(x, xy^2+y^4+y^9). 
$$
This completes the proof. 
\qed

\

\noindent {\sl Proof of 3 and 4 in Proposition \ref{singdeg thm}} :\; 
The proof is similar to that of 2 in Proposition \ref{singdeg thm}.
{ We can directly show that 
$$ 
\textstyle 
(x, xy^2+y^5+by^7+cy^8+dy^9) \sim_{\A^9} (x, xy^2+y^5+(d-\frac{5}{6} b^2)y^9)
$$
by suitable coordinate changes.

On the other hand 
$$\displaystyle
(x,xy^2+y^6+y^7+O(8))\sim_{\A}(x,xy^2+y^6+y^7+\alpha y^9+O(10))
$$ 
for $\alpha\in\R$, is shown by Rieger in the proof of Lemma 3.2.1:3 in \cite{Rieger}. 
This completes the proof. }
\qed

\subsection{{\bf Case 2}:  $d\lambda(0)=0$,    $\rk H_\lambda (0)=1$}

We deal with types $4_k\, (k=3,4,5), 16-18$ of $\Acod = 4, 5, 6$. 

\begin{prop}\label{doubdeg thm}
 For a plane-to-plane map-germ $f$ of corank one, 
\begin{enumerate}
\renewcommand{\labelenumi}{(\arabic{enumi})}
\item
$j^3f(0) \sim_{\A^3} (x,y^3)\Longleftrightarrow$
$$
d \lambda (0) = 0,   \;\; \rk H_\lambda (0) =1, \;\;
 {\eta} ^2 \lambda (0) \neq 0.
$$
\item
If we write 
$f =(x,y^3+\sum _{i+j\geq 4}a_{ij} x^iy^j)$, 
$$\textstyle
\begin{array}{ccl}
\vspace{2mm}
a_{31} \not=0   &\Longleftrightarrow&  f \sim_{\A}(x,y^3+x^3y)\cdots \fbox{\mbox{$4_3$}}\,, \\
\vspace{2mm}
\left\{\begin{array}{l}
a_{31} =0  \\
a_{41}-\frac13 a_{22}^2 \not=0
\end{array}\right.  &\Longleftrightarrow&  f \sim_{\A}(x,y^3\pm x^4y)\cdots \fbox{\mbox{$4_4$}}\,, \\
\left\{\begin{array}{l}
a_{31} =a_{41}-\frac13 a_{22}^2 =0\\
a_{51}-\frac23a_{32}a_{22}+\frac13a_{13}a_{22}^2\not=0
\end{array}\right.  &\Longleftrightarrow&  f \sim_{\A}(x,y^3\pm x^5y)\cdots \fbox{\mbox{$4_5$}}\,.
\end{array}
$$
\item
$j^4f(0) \sim_{\A^4} (x,x^2y+y^4)$ $\Longleftrightarrow$
$$
d \lambda (0) = 0, \;\; \rk H_\lambda (0) =1, \;\;
{\eta} ^2 \lambda (0) =0, \;\; {\eta} ^3 \lambda (0) \neq 0. 
$$
\item
 If we write 
$f =(x,x^2y+y^4+\sum _{i+j\geq 5}a_{ij} x^iy^j)$, 
\begin{eqnarray*}
a_{05} \not=0   &\Longleftrightarrow&  f \sim_{\A}(x,x^2y+y^4\pm y^{5})\cdots \fbox{16}\,, \\
a_{05} =0   &\Longleftrightarrow&  f \sim_{\A}(x,x^2y+y^4)\;\; \;\;\;\;\;\; \cdots \fbox{17}\,. \;\;  
\end{eqnarray*}
\item
 $j^4f(0) \sim_{\A^4} (x,x^2y+xy^3) \;\; \mbox{or} \;\; (x,x^2y) $ 
$\Longleftrightarrow$
$$
d \lambda (0) = 0, \;\; \rk H_\lambda (0) =1, \;\;
{\eta} ^2 \lambda (0) ={\eta} ^3 \lambda (0) = 0.
$$
\item
If we write 
$f =(x,x^2y+xy^3+\sum _{i+j\geq 5}a_{ij} x^iy^j)$, 
$$\textstyle
\left\{\begin{array}{l}
\vspace{1mm}
a_{05}\not= \frac32, \frac{9}{5}  \\
a_{06}(5 a_{05}-9)-15 a_{14} a_{05}  \neq 0 
\end{array}\right.
$$
$$
\Longleftrightarrow  j^7f(0) \sim_{\A^7}(x,x^2y+xy^3+ \alpha y^{5}+y^6+\beta y^7 ),
$$
and excluding exceptional values of $\alpha$ and $\beta$,
$$
f \sim_{\A}(x,x^2y+xy^3+ \alpha y^{5}+y^6+\beta y^7 )\;\;\; 
 \cdots \fbox{18}\,. \\
$$
 \end{enumerate}
\end{prop}

\vspace{1mm}

 Note that  we exclude the type $(x,x^2y)$, because it has codimension $7$, 
while the type $18$ has codimension $6$. 



\

\noindent{\it Proof of 1, 3 and 5 in Proposition \ref{doubdeg thm}.} 
{ We can prove these statements by similar way to  proof of 1 in Proposition \ref{nondeg thm}.
\qed}

\

\noindent{\it Proof of 2 in Proposition \ref{doubdeg thm}.} 
{ Let 
$$\textstyle 
f(x,y) =(x,y^3+\sum _{i+j\geq 4}a_{ij} x^iy^j).$$
By a coordinate change of the source plane $f$ is equivalent to }
$$\textstyle 
(x,y^3+ a_{31}x^3y+\sum_{i+j\ge5}{a'}_{ij}x^iy^j)
$$
for some $a'_{ij}$ where ${a'}_{41}=a_{41}-\frac13 a_{22}^2-\frac13a_{13}a_{31}$.
Since $4_3$-type is $4$-determined, 
$a_{31}\not= 0$ leads to 
$$f \sim_{\A}(x,y^3\pm x^3y).$$

Suppose $a_{31}= 0$. In entirely the same way as above, 
$f$ is equivalent to 
$$
\textstyle 
(x,y^3+ (a_{41}-\frac13 a_{22}^2)x^4y+\sum_{i+j\ge6}{b}_{ij}x^iy^j)
$$
for some $b_{ij}$.
Since $4_4$-type is $5$-determined,  $a_{41}-\frac13 a_{22}^2\not=0$ leads to
$$f \sim_{\A}(x,y^3\pm x^4y).$$

If $a_{41}-\frac13 a_{22}^2=0$, then 
$f$ is equivalent to 
$$
\textstyle 
(x,y^3+ (a_{51}-\frac23a_{32}a_{22}+\frac13a_{13}a_{22}^2)x^5y+O(7)). 
$$
Since $4_5$-type is $6$-determined, the claim 2 follows. \qed

\

\noindent{\it Proof of  4 in Proposition \ref{doubdeg thm}.} 
Let 
$$\textstyle 
f(x,y) =(x,x^2y+y^4+\sum _{i+j\geq 5}a_{ij} x^iy^j).$$
Rewrite variables, and substitute 
\begin{eqnarray*}
&&\tilde{x}=x\\
&&\textstyle
\tilde{y}=y- \sum_{i+j=5,\; i\ge2 }a_{ij}{x}^{i-2}{y}^j,
\end{eqnarray*}
then we see that $f$ is equivalent to 
$$\textstyle 
(x,x^2y+y^4+a_{14}xy^4+a_{05}y^5+O(6)). 
$$ 
Now we show 
$$ 
(x, x^2y+y^4+cxy^4+dy^5) \sim_{\A^5}(x, x^2y+y^4+dy^5).
$$
This is explicitly given by 
$\tilde{x}=x$ and 
$\tilde{y}=y - \frac {c } {3}x y$ 
and the coordinate change of the target:
$$\textstyle 
(X,Y)\mapsto \left(X, Y + \frac {c}{3}X Y + \frac { c^2} {9} X^2 Y + \frac {c^3} {27} 
 X^3 Y \right).
$$
Since 16 and 17-types are $5$-determined, 
the claim is proved.
\qed

\

\noindent{\it Proof of  6 in Proposition \ref{doubdeg thm}.} 
At first, for $d\not=\frac95$
\begin{eqnarray*}
&&(x,x^2y+xy^3+cxy^4+dy^5+exy^5+gy^6)\\
&&\qquad 
\textstyle 
\sim_{\A^6}
(x,x^2y+xy^3+dy^5+Pxy^5+(g-\frac{15cd}{5d-9})y^6)
\end{eqnarray*}
holds where $P$ is a constant. { This is also given by some complicated coordinate changes.} 

Next we see 
$$
(x,x^2y+xy^3+dy^5+exy^5+gy^6)
\sim_{\A^6}
(x,x^2y+xy^3+dy^5+gy^6)
$$ 
for $d\not=\frac32$.
This follows from 
\begin{eqnarray*}
&&\tilde{x}=x\\
&&\textstyle 
\tilde{y}=y + \frac {e } {2 (3 - 2 d)} x y+ \frac {e^2 } {4 (3 - 2 d)^2} x^2 y+ \frac {e^3 } {8 (3 - 
2 d)^3}x^3 y - \frac {e } {3 - 2 d} y^3- \frac {5 e^2} {4 (3 - 2 d)^2} x y^3
\end{eqnarray*}
and $\textstyle
(X,Y)\mapsto \left(X, Y + \frac { e} {2(2d-3)}XY \right)$. 
Now let 
$$\textstyle 
f(x,y)=(x,x^2y+xy^3+\sum _{i+j\ge5}a_{ij} x^iy^j).$$
{ A simple coordinate change shows 
that $f$ is equivalent to}  
$$\textstyle 
(x,x^2y+xy^3+a_{14}xy^4+ a_{05}y^5+Qxy^5+a_{06}y^6+O(7))
$$
where $Q$ is a constant, and hence by the above argument,  
$f$ is equivalent to 
$$\textstyle 
(x,x^2y+xy^3+ a_{05}y^5+(a_{06}-\frac{15 a_{14} a_{05}}{5 a_{05}-9}) y^6+O(7)) 
$$
for $a_{05}\not= \frac32, \frac{9}{5}$.
 Finally by a similar coordinate change as above again, $f$ is equivalent to 
$$\textstyle 
(x,x^2y+xy^3+ a_{05}y^5+(a_{06}-\frac{15 a_{14} a_{05}}{5 a_{05}-9}) y^6+b_{16}xy^6+b_{07}y^7+O(8)) 
$$
for some $b_{ij}$; 
Then $xy^6$ in the second component is killed, 
and moreover, if the coefficient of $y^6$ is not zero, then 
$f$ is $\A$-equivalent to 
$$\textstyle 
(x,x^2y+xy^3+ \alpha y^5+y^6+\beta y^7+O(8)) 
$$
that follows from the same argument as in the proof of Proposition 3.2.2:2 in \cite{Rieger}. 
For generic values $\alpha$ and $\beta$, 
it implies the claim 6, by the determinacy result. 
This completes the proof. \qed

\begin{rem}\label{4_k}
{ \rm $4_k$-type can be characterized  as follows:  For a plane-to-plane map-germ $f$ of corank one, 
\begin{center}
$f \sim_{\A} (x, y^3\pm x^ky) $
$\Longleftrightarrow$ $\lambda$ is $A_{k-1}$-type and $\eta^2 \lambda(0)\neq 0$.
\end{center}
This special feature of $4_k$-type would be explained with the augmentation-theory (see \cite{ORW1}). This will be studied somewhere else.  
}
\end{rem}

\subsection{{\bf Case 3}:  $d\lambda(0)=0$,    $\rk H_\lambda (0)=0$}
 Finally we deal with type $19$ of $\Acod =6$.
\begin{prop}\label{trideg them}
 For a plane-to-plane map-germ $f$ of corank one, 
\begin{enumerate}
\renewcommand{\labelenumi}{(\arabic{enumi})}

\item
$j^4f(0) \sim_{\A^4} (x, x^3y+\alpha x^2y^2+y^4) \; \mbox{or} \; (x, \alpha x^2y^2+y^4)$ 
$\Longleftrightarrow$
$$
d \lambda (0) = 0, \;\;  \rk \, H_\lambda (0) =0, \;\; \eta^3\lambda(0) \neq 0. 
$$
\item {\rm (Rieger \cite{Rieger}) }
 If we write 
$f =(x, x^3y+\alpha x^2y^2+y^4+\sum _{i+j\geq 5}a_{ij} x^iy^j)$, 
$$
\Delta = 8 \alpha a_{41}-12a_{32}-4\alpha^2a_{23}+4\alpha a_{14}+(3+2\alpha^3)a_{05}
\not=0  
$$
 $$
 \Longleftrightarrow  f \sim_{\A}(x, x^3y+\alpha x^2y^2+y^4+x^3y^2)\cdots \fbox{19} 
$$
 \end{enumerate}
\end{prop}

 Note that the $\A$-codimension of $(x, \alpha x^2y^2+y^4)$ is greater than $6$, 
so we exclude it, while the type $19$ has codimension $6$.  {The claim 2 is due to Rieger, that can be seen in the proof of Prop. 3.2.3.1 in \cite{Rieger}, and we can show the first claim in the same manner with the previous sections. }

\


\section{Application to projection of surface in $3$-space} 
This section is devoted to an application of our criteria 
to singularities arising in parallel/central projection of a surface in $3$-space. 
Our standing point is to look at this problem as a typical one of 
{\it $\A$-recognition problem of plane-to-plane map-germs 
arising in a concrete geometric setting}. 

\subsection{Parallel and central projections}\label{Parallel-and-central-projections}
Let $\iota: M \hookrightarrow \R^3$ be an embedding of smooth surface. 
\begin{defin}
 {\it A parallel projection} of a smooth surface $M$ to the plane 
 is the restriction to $M$ of a linear orthogonal projection $pr: \R^3\rightarrow \R^2$.  
 \end{defin}

The direction of orthogonal projection has two dimensional freedom; 
the space of directions is just the $2$-dimensional sphere $S^2$. 
There is naturally produced 
a $2$-parameter family of parallel projections, $M\times U \to \R^2$, 
where $U$ is any small open subset of $S^2$.  
Hence a na\"Ive guess is that 
any plane-to-plane germs of $\A$-$\mbox{cod}\, \leq 4$ might appear generically 
in parallel projection of surface $M$ at some points. 
In fact it is true. 
\begin{thm}[{ Arnold \cite{Arnold}, Gaffney-Ruas \cite{Gaffney-Ruas, Gaffney}, Bruce \cite{Bruce} }] \label{bruce} 
For a generic surface $M$, singularities arising in parallel projections of $M$  are $\A$-equivalent to the germs of $\A$-$\mbox{cod}\, \leq 4$  in  Table \ref{Riegerlist}. 
\end{thm}

\begin{rem}{\rm 
We should  remark about what  the word ``generic"  means. 
A precise statement is as follows:  
there exists a residual subset of the space of all embeddings of $M$ into $\R^3$ 
 (equipped with $C^\infty$-topology) so that 
 for each element $\iota: M \hookrightarrow \R^3$ of this subset, 
 any parallel projection $pr|_M: M \to \R^2$ admits only singularities 
 of  $\A$-$\mbox{cod}\, \leq 4$ listed  in  Table \ref{Riegerlist}. 
 Below we abuse  the word ``generic" in the same manner for several similar situations; 
Perhaps that would not cause any confusion. 
}
\end{rem}

\begin{defin}
A {\it central projection} of $M$ from a viewpoint $p=(a,b,c) \in \R^3-M$ is defined 
by the restriction 
$$\varphi_p:=\pi_p|_M: M \to \R P^2$$
of the canonical surjection on the projective plane  
$$\pi_p: \R^3 -\{p\}\to \R P^2, \quad x \mapsto \mbox{line generated by $x-p$}.$$
\end{defin}

There is $3$-dimensional freedom of the choice of viewpoint $p$; 
there is naturally produced 
a $3$-parameter family of central projection, $M \times U \to \R P^2$, 
where $U$ is any small open subset of the complement $\R^3-M$.  
Therefore we might have expected  
that any plane-to-plane germs of $\A$-$\mbox{cod}\, \leq 5$
would appear in central projection generically. 
However it is not the case. 
Arnold and Platonova proved the following remarkable theorem \cite{Arnold, Platonova}:

\begin{thm}[Arnold \cite{Arnold}, Platonova \cite{Platonova}]\label{APthm}
For a generic surface $M$, and for any $p \in \R^3$ not lying on $M$, 
the germ $\varphi_p: M, x \to \R P^2, \varphi_p(x)$ at any point $x \in M$ is 
$\A$-equivalent to one of the list of germs with $\A$-$\mbox{cod}\, \leq 5$
in Table \ref{Riegerlist} 
except for $12, 16$ and unimodal type $8$. 
\end{thm}

So the three types $12, 16$ and $8$ are excluded in the list  
of singularities arising in central projection of a generic surface, 
in other words, this geometric setting makes a strong restriction 
on the appearance of singularities of plane-to-plane germs of $\A$-$\mbox{cod}=5$. 
Our  criteria are applied to detecting $\A$-types of map-germs 
arising in this special geometric setting. 
Then we give not only a new transparent proof of Theorem \ref{APthm} 
 in the context of Rieger's classification but also 
some extension as stated in the following theorem: 

\begin{thm}\label{mainthm2} 
For a generic one-parameter family of embeddings $M \times I \to \R^3$, $(x, t) \mapsto \iota_t(x)$, 
the central projection $\pi_p\circ \iota_t: M \to \R P^2$  for any $t$ and any viewpoint $p$ 
admit only $\A$-types with $\Acod \le 5$ and 
types $12, 16, 8, 4_5, 9, 11_9, 13, 17, 19$ with $\Acod \, 6$. 
Namely, each type of $10,15,18$ with $\Acod \, 6$ does not appear generically. 
\end{thm}

In Rieger \cite{Rieger3}, parallel projection of moving surfaces with one-parameter 
has been considered. Theorem \ref{mainthm2} generalizes it in a much more general form.

\subsection{Proofs of Theorems \ref{APthm} and \ref{mainthm2}} 
First we explain the main idea of the proof: This is a slightly modified version of the method 
 which W. Bruce used for the proof of Theorem \ref{bruce} in  \cite{Bruce}. 
Now our setting is for central projection. 
Let $(x,y,z)$ be the coordinates in $\R^3$. 
Assume that 
$x_0=(1,0,0) \in M$ and  $T_{x_0}M$ is $xy$-plane, 
and viewpoints   
$p=(a,b,c) \in \R^3$ is close to $0$ enough. 

Then $M$ is locally expressed by its 
{\it Monge form} $z=f(x,y)$ centered at  $x_0 \in M$, i.e., 
$$M =\{(1+x, y, f(x,y)) \in \R^3\}$$
as a set-germ at $x_0$ and $df(0,0)=0$. 
Note that infinitesimal information of $M$ at $x_0$ can be deduced from 
the Taylor expansion of $f$.
In particular, the germ at $x_0$ 
of the central projection from the viewpoint $p$, 
$$\varphi_{p,f}: M, x_0 \to \R P^2, y_0,$$
($y_0=\varphi_{p,f}(y_0)$) is explicitly written by 
$$
\varphi_{p,f}(x,y):=\left(\frac{y-b}{1+x-a},\frac{f(x,y)-c}{1+x-a}\right)
$$
using local coordinates $(x,y)$ of $M$ and $[1:X:Y]$ of $\R P^2$. 

Let $V_\ell$  be the $\ell$-jet space of the Monge form $z=f(x,y)$ at $0$ 
 (that is the space of polynomials of degree greater than $1$ and less than or equal to $\ell$). 
Also denote by $J^\ell (2,2)$ the jet space of $\R^2,0 \to \R^2,0$. 
We then define $\Phi (p, j^\ell f(0)) \in J^\ell (2,2)$ to be the $\ell$-jet of 
$\varphi_{p,f}(x,y)-\varphi_{p,f}(0,0)$ at the origin, and consider the following diagram:
$$
\xymatrix{
\R^3 \times V_\ell \ar[r]^{\Phi} \ar[d]_{pr} & J^\ell(2,2), & (p, j^\ell f(0)) \ar[r]^\Phi \ar[d]_{pr} & j^\ell \varphi_{p, f}(0). \\
V_\ell & & j^\ell f(0)
}
$$

Note that $J^\ell (2,2)$ is stratified by $\A^\ell$-orbits 
(those strata of low codimension  are given in Rieger's list). 
Therefore $\Phi$ induces a stratification of $\R^3 \times V_\ell$. 

\begin{defin}
For an $\A^\ell$-orbit  $W \subset J^\ell(2,2)$, define 
$$G_{W}:= pr(\Phi^{-1}({W})) \;\; \subset \;\;  V_\ell.$$
\end{defin}

Since any $\A^\ell$-orbit $W$ is a semi-algebraic subset of $J^\ell(2,2)$, 
 $\Phi^{-1}({W})$ and hence $G_{W}$ turns out to be semi-algebraic.

Next we discuss a certain variant of Thom's transversality theorem. 
First we fix an Euclidean metric of $\R^3$ and its orientation. 
Suppose that we are given an embedding $M \subset \R^3$ with the unit normal vector field 
$\bn:M \to \R^3$. 
Let $U$ be an  open subset of $M$ with 
unit tangent vector fields $\bv$ and $\bw$ so that 
$\{\bv,\bw,\bn\}$ is an orthonormal frame with respect to the fixed orientation. 
At each point $q \in U$,  linear coordinates $x_q, y_q, z_q$ of $\R^3$ centered at $q$ 
are chosen 
to coincide with the oriented lines defined by vectors $\bv(q),\bw(q),\bn(q)$, respectively.  
Writing $U$ near $q$ as the Monge form $z_q=f_q(x_q, y_q)$, 
we associate to each point $q$ the Taylor expansion of $f_q$ truncated to degree $\ell$: 
Then 
$$\Theta: U \to V_\ell, \quad \Theta(q):=j^\ell f_q(0)$$ 
is defined 
after rewriting the variables such as $x_q=x$, $y_q=y$, $z_q=z$.

Globally, we take an open cover $\{U_i\}$ of $M$ so that for each subscript $i$, 
we have $\Theta_i: U_i \to V_\ell$  in the same way just as mentioned. 
Note that there is the right linear action on $V_\ell$ of the rotation group $SO(2)$: 
For $q \in U_i \cap U_j$, 
the difference between $\Theta_i$ and $\Theta_j$ at $q$ is only caused by 
this action. Then, 
in \cite[Thm.1]{Bruce}, 
the following version of transversality theorem is proved:

\begin{prop}[Bruce \cite{Bruce}]
Let $X\subset V_{\ell}$ be an $SO(2)$ invariant submanifold. 
For generic surface $M$ in $\R^3$, $ \Theta_i: U_i \to V_\ell$ is transverse to $X$.
\end{prop}

Obviously, $G_W$ is an $SO(2)$-invariant subset. 
It immediately implies the following assertion by a standard argument of transversality theorem:  

\begin{cor}\label{projcor}
(1) If $\mbox{codim}\, G_{W} \ge 3$, then 
 for a generic embedded surface $M$, the central projection $\varphi_p: M \to \R P^2$ 
from any viewpoint $p$ does not admit ${W}$-type singularity at any point of $M$.  
(2) For a generic $s$-parameter family of embeddings of $M$ into $\R^3$, 
any central projection admits only singularities of type $W$  
with $\mbox{codim}\, G_{W} \le 2+s$. 
\end{cor}

From Corollary \ref{projcor}, 
our main task for proving Theorem \ref{mainthm2} 
is to determine $\mbox{codim}\, G_W$ for all $W$ in consideration. 
To do this,  we describe explicitly the defining equations of $G_W$.
We obtain the following result:

\begin{prop}\label{projthm}
 Table \ref{List4} is the list of $\codim G_W$ for all the map-germs of $\A$-$\codim\le6$, with $\ell$ large enough. 
 In addition, $\codim G_W\ge4$ holds for all the map-germs of $\Acod\ge7$. 
\end{prop}

Theorems \ref{APthm} and 
 \ref{mainthm2} immediately follow from Proposition \ref{projthm} and  Corollary \ref{projcor}. 

\begin{table}
$$
\begin{array}{c|c| l  }
\mbox{cod}\, G_W & \A\mbox{-cod}  &\mbox{type}\,  W \\
\hline
\hline
0&0&1 \\
&1&2  \\
&2&3  \\
\hline 
1&3&4_2 , 5  \\
&4&4_3 \\
\hline 
2&4&6, 11_5 \\
&5&4_4, 7, 11_7  \\
\hline 
3&5&12,16, 8 \\
&6&4_5, 9, 11_9,13,17,19\\
\hline
4&6&10,15,18 \\
\end{array}
$$
\caption{Codimension of $G_W$ and $\A$-codimension of $W$}
\label{List4}
\end{table}

To show Proposition \ref{projthm},  
we consider the same cases as in the previous section: cases $0, 1, 2, 3$. 
The proof will be done as follows. 
From now on we write 
$$ f(x,y)=\sum_{ i+j \ge 2} c_{ij} x^i y^j.$$
For each $\A$-type in Table \ref{Riegerlist}, 
we will apply our criteria in Chapter 2 to the plane-to-plane germ of the following form 
$$
\varphi_{p,f}(x,y)=
\left(\frac{y-b}{1+x-a},\frac{\sum  c_{ij} x^i y^j-c}{1+x-a}\right). 
$$
Then we obtain a certain condition in variables 
$$a, b, c, c_{20}, c_{11}, c_{02}, c_{30}, c_{21}, \cdots $$ 
so that $\varphi_{p,f}$ is $\A$-equivalent to the $\A$-type. 
That is nothing but the condition defining 
the semi-algebraic subset $\Phi^{-1}(W)$ 
in $\R^3 \times V_\ell$ 
for the corresponding $\A^\ell$-orbit 
$W \subset J^\ell(2,2)$ (with $\ell$ larger than 
the determinacy order). 
The condition consists of polynomial equations 
and  inequalities. Simply we call the (system of) equations 
the {\it defining equation  of $\Phi^{-1}(W)$}. 
By eliminating the variables $a, b, c$ from the  equation, 
we obtain the defining equation of  $G_W$. 
The inequalities do not affect the codimension. 

In general the codimension of $\Phi^{-1}(W)$ is equal to that of $W$, 
therefore the main task is to check how the projection 
$pr$ affects the defining equation of $G_W$.

\subsection{\bf Case 0}
Here we think of  the case $d\lambda \not= 0$. 
This { case} automatically implies $\eta \lambda(0)=0$ in common other than 2-type.   
By this condition, it follows that
$$
(1-a)^2 c_{20} +(1-a)(-b) c_{11}+(-b)^2 c_{02}=0
$$
and $\det \mbox{Hess}f(0) <0$,  that is, $f$ is hyperbolic at  the origin. 
Hence from now on we assume that 
$$c_{20}=c_{02}=0, \;\;\mbox{and} \quad c_{11}\not=0 $$ 
by taking a suitable coordinate $x, y$ via $SO(2)$-action. 
With this condition, our calculations become much easier.

Let us determine the defining equation of $\Phi^{-1}(W)$, $G_W$  and $\codim \, G_W$ for each $\A^\ell$-orbit $W$. 

\

\noindent
\fbox{$2$} \; 
$\Phi^{-1}(W)$ is given by 
$c=0$. 
Therefore there is no defining equation for $G_{W}$.
Hence  $\codim \, G_{W}=0$.

\

\noindent
\fbox{$3$} \;  
 $\Phi^{-1}(W)$ is given by 
$b=c=0$. 
Therefore there is no defining equation for $G_{W}$.
Hence  $\codim \, G_{W}=0$.

\

\noindent
\fbox{$5$} \; 
 $\Phi^{-1}(W)$ is given by 
$b=c=0$ and $c_{30}=0$. 
 Then $G_{W}$ is given by 
$c_{30}=0$. 
Hence  $\codim \, G_{W}=1$.

\

\noindent
\fbox{$6$} \; 
$\Phi^{-1}(W)$ is given by 
$b=c=0$ and $c_{30}=c_{40}=0$. 
Then $G_{W}$ is given by 
$c_{30}=c_{40}=0$. 
Hence  $\codim \, G_{W}=2$.

\

\noindent
\fbox{$7$} \;  
$\Phi^{-1}(W)$ is given by 
$b=c=0$, 
$c_{30}=c_{40}=0$ and 
$A(1-a)^2+B(1-a)+C=0$
where $A, B, C$ are some polynomials in $c_{ij}$. 
We get this equation from  the additional conditions of Taylor expansions: $a_{07}-\frac58 a_{06}^2=0$ in our criterion 
for $7$-type 
(see Proposition \ref{nondeg thm}):  Under the condition
$$
d\lambda(0)\not=0, \eta\lambda(0)=\eta^2\lambda(0)=\eta^3\lambda(0)=0, \eta^4\lambda(0)\not=0,
$$
$\varphi_{p,f}$ is $\A$-equivalent to $(x, xy+y^5 +\sum_{i+j\ge6} a_{ij}x^iy^j) $
where
$$
a_{06}=\frac{1}{c_{11}c_{50}(1-a)} 
\left\{
(-5c_{21}c_{50}+6c_{11}c_{60})(1-a)-c_{11}c_{50}
\right\},
$$
and
 \[
a_{07}=\frac{1}{c_{11}^2c_{50}(1-a)^2}
\left\{
{\small \begin{array} {l}
(20c_{21}^2 c_{50}-5c_{11}c_{31}c_{50}-6c_{11}c_{21}c_{60}+c_{11}^2c_{70})(1-a)^2\\
 +(6c_{11}c_{21}c_{50}-c_{11}^2c_{60})(1-a)+c_{11}^2c_{50}
\end{array}}
\right\}.
\]
(See  Proof of 1 in Proposition \ref{nondeg thm}.)
Here $a_{07}-\frac58 a_{06}^2=0$ gives an equation in the variable $a$ 
with $A, B, C$ depending only on $c_{ij}$: 
$$A(1-a)^2+B(1-a)+C=0.$$
In addition, 
$A$ (also $B$, $C$) is independent from the other two equations $c_{30}=c_{40}=0$ ; 
for instance, we see that there is a monomial $c_{11}c_{31}c_{50}^2$ in $A$. 
The variable $a$ is solved in $c_{ij}$ generically; 
the locus in $V_\ell$ where $a$ is not solved is defined by $A=B=C=0$, 
but it has high codimension, 
so this quadratic equation does not affect $\codim \, G_{W}$.
Thus, $G_{W}$ is given by 
$c_{30}=c_{40}=0$. 
Hence  $\codim \, G_{W}=2$. 

\
 
\begin{rem}\label{h-focal}{\rm 
For a generic surface,  hyperbolic points where 
the Monge form satisfies that $c_{30}=c_{40}=0$ are isolated, 
since $G_W$ has codimension $2$ in $V_\ell$.  
Look at such a point of the surface from a viewpoint lying on the $a$-axis ($b=c=0$), 
then 
the central projection produces the butterfly singularity (6-type). 
However,  there is an exception: from at most two points on the $a$-axis 
which are given by the solution $a=a(c_{ij})$ of 
the quadric equation, the central projection admits the elder-butterfly singularity (7-type). 
 These exceptional points are called {\it $h$-focal points} (``$h$'' for ``hyperbolic'')  by Platonova \cite{Platonova}}.
\end{rem}

\

\noindent
\fbox{$8$} \; 
$\Phi^{-1}(W)$ is given by 
$b=c=0$ and $c_{30}=c_{40}=c_{50}=0$. 
Then $G_{W}$ is given by 
$c_{30}=c_{40}=c_{50}=0$. 
Hence  $\codim \, G_{W}=3$.
The difference between 7-type and 8-type (although they have the same $\A$-codimension) 
  is the difference between closed conditions $a_{07}-\frac58 a_{06}^2=0$ and $\eta^4 \lambda(0)=0$. 
As mentioned in Remark \ref{h-focal}, the former condition on coefficients   
determines the position of viewpoint,  
while the geometric condition $\eta^4\lambda(0)=0$ is that of the Monge form. 
Also in the following other calculations, 
this kind of difference makes the difference of $\codim \, G_W$.

\

\noindent
\fbox{$9$} \; 
$\Phi^{-1}(W)$ is given by 
$b=c=0$, $c_{30}=c_{40}=c_{50}=0$, 
$A(1-a)^2+B(1-a)+C=0$
for some polynomials $A, B, C$ in $c_{ij}$. 
{ The last equation} comes from the condition $a_{08}-\frac{3}{5} a_{07}^2=0$ in our criterion (Proposition \ref{nondeg thm}) in the similar way as the case of no. $7$ 
 (e.g., $A$ contains the monomial $c_{11}c_{31}c_{60}^2$, so it is independent from other three equations). 
Then $G_{W}$ is given by 
$c_{30}=c_{40}=c_{50}=0$. 
Hence  $\codim \, G_{W}=3$.

\

\noindent
\fbox{$10$}\: 
$\Phi^{-1}(W)$ is given by 
$b=c=0$, $c_{30}=c_{40}=c_{50}=c_{60}=0.$
Then $G_{W}$ is given by 
$c_{30}=c_{40}=c_{50}=c_{60}=0$. 
Hence  $\codim \, G_{W}=4$.

\subsection{\bf Case 1}
Here we think of  the case that 
$d\lambda (0)= 0$ and $H_\lambda (0)$ is non-degenerate. 
Note that we can always assume that 
$c_{11}=0$ by taking a suitable rotation of $xy$-plane. 
Then the condition  $\frac{\rd}{\rd x}\lambda (0)=\frac{\rd}{\rd x}\lambda (0)=0$ leads to 
$$
(1-a)c_{20}=0, \quad b\;c_{02}=0, 
$$
i.e., $c_{20}=b=0$ or $c_{20}=c_{02}=0$, 
hence $f$ is parabolic or umbilic (flat) at the origin, 
while $\rk \, H_\lambda (0)=2$ leads to an open condition. 
Note that the locus consisting of umbilic Monge forms is defined 
by vanishing the $2$-jet, so it has codimension $3$  in $V_\ell$. 
Below, for the simplicity, we treat only with the case that $c_{20}=b=0$ 
(i.e., parabolic forms); 
In fact, taking $c_{02}=0$ instead of $b=0$ (i.e., considering umbilic forms) 
does not affect  the codimension of $G_W$ as seen in $11_5$-type.

\

\noindent
\fbox{$ 4_2$} \; 
$\Phi^{-1}(W)$ is given by 
$b=c=0$ and $c_{20}=0$. 
Then $G_{W}$ is given by 
$c_{20}=0$. 
Hence  $\codim \, G_{W}=1$. 
 That is, the locus of parabolic Monge forms has codimension $1$  in $V_\ell$. 

\

\noindent
\fbox{${11}_5$}\; 
$\Phi^{-1}(W)$ is given by 
$b=c=0$ and $c_{20}=c_{30}=0$. 
Then $G_{W}$ is given by 
$c_{20}=c_{30}=0$. 
Hence  $\codim \, G_{W}=2$. 

 Precisely saying, the component of $G_W$ having parabolic Monge forms 
has codimension $2$ in $V_\ell$, where  $\eta^2 \lambda(0)=0$ leads to 
the extra equation $c_{30}=0$. 
If we take $c_{02}=0$ instead of $b=0$, i.e., we consider the umbilic Monge form, 
then $\eta^2 \lambda(0)=0$ implies an equation of $a$, $b$ and $c_{ij} \, (i+j=3)$, 
instead of  $c_{30}=0$. 
So we can generically solve  $a$ or $b$ in $c_{ij}$. 
Thus the component of $G_W$ having umbilic Monge forms 
remains to be of codimension $3$ in $V_\ell$ 
and it does not affect the codimension of $G_W$. 
This argument is also valid in the following other types, so we will not repeat it below.

\

\noindent
\fbox{${11}_7$}\; 
$\Phi^{-1}(W)$ is given by 
$b=c=0$, $c_{20}=c_{30}=0$, $A(1-a)+B=0$
for some polynomials $A, B$ in $c_{ij}$. 
Then $G_{W}$ is given by 
$c_{20}=c_{30}=0$. 
Hence  $\codim \, G_{W}=2$.

Remark that the equation $A(1-a)+B=0$ 
arises from the condition $a_{05}=0$  in our criterion (see Proposition \ref{singdeg thm}) 
after rewriting $\varphi_{p, f}$ to be $(x, xy^2+\sum_{i+j\ge 5} a_{ij}x^iy^j)$ 
by an explicit coordinate change. 

\begin{rem}\label{p-focal}{\rm 
As seen in Remark \ref{h-focal}, we also have an exceptional point here.
We look at parabolic points on the surface where $c_{20}=c_{30}=0$ from a viewpoint lying on the $a$-axis ($b=c=0$), that is the unique asymptotic line. Then the gulls singularity ($11_5$-type) appears on the line except for the point $(a,0,0)$ 
where $a$ is given by $A(1-a)+B=0$. This exceptional point is called {\it $p$-focal} point (``$p$'' for parabolic) by Platonova \cite{Platonova}, and at this point 
the ugly-gulls singularity ($11_7$-type) appears.}
\end{rem}

\

\noindent
\fbox{${11}_9$}\; 
\noindent
$\Phi^{-1}(W)$ is given by 
$b=c=0$, $c_{20}=c_{30}=0$, $A(1-a)+B=0$, 
$a_{07}-2a_{15}+4 a_{23}=0$
 where $a_{ij}$ are functions in $c_{ij}$'s and $a$ obtained  
in entirely the same way as the above case $11_7$. 
Solve the variable $a$ by $A(1-a)+B=0$, 
and then the last equation yields a non-trivial equation, say $ C(c_{ij})=0$, 
which is independent from other equations. 
Therefore $G_{W}$ is given by 
$c_{20}=c_{30}=0$ and $C(c_{ij})=0$. 
Hence  $\codim \, G_{W}=3$.

\

\noindent
\fbox{$12$}\; 
$\Phi^{-1}(W)$ is given by 
$b=c=0, \;c_{20}=c_{30}=c_{40}=0$. 
Then $G_{W}$ is given by 
$c_{20}=c_{30}=c_{40}=0$. 
Hence  $\codim \, G_{W}=3$. 

\

\noindent
\fbox{$13$}\; 
$\Phi^{-1}(W)$ is given by 
$b=c=0$, $c_{20}=c_{30}=c_{40}=0$, $A(1-a)+B=0$
for some polynomials $A, B$ in $c_{ij}$. 
Here the last equation comes from the condition $a_{06}=0$ 
in our criterion (see Proposition \ref{singdeg thm}). 
Then $G_{W}$ is given by 
$c_{20}=c_{30}=c_{40}=0$. 
Hence  $\codim \, G_{W}=3$.

\

\noindent
\fbox{$15$}\; 
$\Phi^{-1}(W)$ is given by 
$b=c=0$, $c_{20}=c_{30}=c_{40}=c_{50}=0$. 
Then $G_{W}$ is given by 
$c_{20}=c_{30}=c_{40}=c_{50}=0$. 
Hence  $\codim \, G_{W}=4$.

\subsection{\bf Case 2}
Here we think of the case $d\lambda (0)= 0$ 
and $H_\lambda (0)$ is degenerate ($\rank=1$). 
Since $d\lambda (0)= 0$, as seen in  case 1, 
we can assume $c_{11}=0$. 
Remark again that in this condition, $d \lambda (0) =0$ leads to $c_{20}=b=0$ 
and 
\[
H_\lambda (0)  = 
-\frac{1}{(1-a)^3}
\left(
    \begin{array}{ccc}
      3(1-a) c_{30}&(1-a)c_{21} \\
      (1-a) c_{21} &(1-a) c_{12} +c_{02} \\
    \end{array}
  \right).
\]
Then $\det H_\lambda (0) =0$ leads to
$$
(3c_{30}c_{12}-c^2_{21})(1-a)+3c_{02}c_{30}=0.
$$
Write it by $C(1-a)+D=0$. 

\

\noindent
\fbox{${4}_3$}\; 
$\Phi^{-1}(W)$ is given by 
$b=c=0$, $c_{20}=0$, $C(1-a)+D=0$.
Then $G_{W}$ is given by 
$c_{20}=0$. 
Hence $\codim G_{W}=1$. 
Here  $C(1-a)+D=0$ does not affect $\codim \, G_{W}$ because 
the condition $C=0$ increases the codimension.

\

\begin{rem}{\rm 
The lips and beaks singularities arise on the unique asymptotic line at a parabolic point;  
there is one exceptional point on the line, given by $C(1-a)+D=0$, where 
the goose singularity appears. This point divides the line into two half lines, 
each of which corresponds to viewpoints for 
either the lips singularity or the beaks singularity of the projection. 
}
\end{rem}

\

\noindent
\fbox{${4}_4$}\; 
$\Phi^{-1}(W)$ is given by 
$b=c=0$, $c_{20}=0$, $C(1-a)+D=0$, $a_{31}=0.$
Here $a_{31}$ is a function in $c_{ij}$ obtained by rewriting 
$\varphi_{p. f}$ to be the suitable form as in Proposition \ref{doubdeg thm}. 
Solve $a$ by $C(1-a)+D=0$ and substitute it into $a_{31}$, 
then denote the resulting equation by $A(c_{ij})=0$, that is independent from others. 
Then $G_{W}$ is given by 
$c_{20}=0$ and $A(a_{ij})=0$. 
 Hence  $\codim \, G_{W}=2$.

\

\noindent
\fbox{${4}_5$}\; 
$\Phi^{-1}(W)$ is given by 
$b=c=0$, $c_{20}=0$,  $C(1-a)+D=0$, $\textstyle a_{31} =a_{41}-\frac13 a_{22}^2 =0.$
Here $a_{ij}$ are obtained by Proposition \ref{doubdeg thm}. 
Then $G_{W}$ is given by 
$c_{20}=0$ and $A(c_{ij})=B(c_{ij})=0$. 
Hence  $\codim \, G_{W}=3$.
$A(c_{ij})$ is the equation $a_{31}=0$, 
and $B(c_{ij})=0$ is the equation $a_{41}-\frac13 a_{22}^2 =0$ after the substitution of $a$.

\

\noindent
\fbox{$16$}\; 
$\Phi^{-1}(W)$ is given by 
$b=c=0, \;c_{20}=c_{30}=c_{21}=0$. 
Then $G_{W}$ is given by 
$c_{20}=c_{30}=c_{21}=0$. 
 Hence   $\codim \, G_{W}=3$.
$c_{30}=c_{21}=0$ come from $\eta^2 \lambda (0) =0$ and $\det H_\lambda (0) =0$. 

\

\noindent
\fbox{$17$}\; 
$\Phi^{-1}(W)$ is given by 
$b=c=0$, $c_{20}=c_{30}=c_{21}=0$, $A(1-a)^2+B(1-a)+C=0$. 
Then $G_{W}$ is given by 
$c_{20}=c_{30}=c_{21}=0$. 
 Hence   $\codim \, G_{W}=3$.
Here $A, B, C$ come from the condition $a_{05}=0$ (see Proposition \ref{doubdeg thm}). 
Hence  this quadratic equation does not affect $\codim \, G_{W}$.

\

\noindent
\fbox{$18$}\; 
$\Phi^{-1}(W)$ is given by 
$b=c=0, \;c_{20}=c_{30}=c_{21}=c_{40}=0$. 
Then $G_{W}$ is given by 
$c_{20}=c_{30}=c_{21}=c_{40}=0$. 
Hence   $\codim \, G_{W}=4$.

\subsection{\bf Case 3}
Here we think of  the case $d\lambda (0)= 0$ and $H_\lambda (0)=O$. In fact, the germ of $\A$-codim$=6$ in (case 3)-singularities is just 19-type.
As seen in the previous sections, we can also suppose $c_{11}=0$.  

\

\noindent
\fbox{$19$}\; 
$\Phi^{-1}(W)$ is given by 
$b=c=0$, 
$c_{20}=c_{30}=c_{21}=0$ and 
$(1-a)c_{12}+c_{02}=0$.
Then $G_{W}$ is given by 
$c_{20}=c_{30}=c_{21}=0$. 
Hence  $\codim \, G_{W}=3$.

\subsection{\bf Types of $\A$-codimension $\ge 7$.} 
We prove the second claim in Proposition \ref{projthm}. 
In Rieger's recognition trees, terminating lines lead to 
$\A$-orbits of $\Acod\ge7$. 
All jets in $J^\ell(2,2)$ indicated by terminating lines satisfy 
closed conditions obtained by replacing 
some inequalities in our criteria by equalities.  
For instance, let us look at 
an orbit $W$ with $\Acod\ge7$ 
over $5$-jets of $\A$-orbit of $(x,y^3)$. 
Any jet $(x,y^3+\sum _{\ell \geq i+j\geq 4}a_{ij} x^iy^j)$ belonging to $W$ satisfies 
the following closed conditions obtained from 
the criterion (2) in Proposition \ref{doubdeg thm} :
\begin{eqnarray*}    
&&d \lambda (0) = 0,   \;\; \rk H_\lambda (0) =1,\\
&&\textstyle a_{31} =a_{41}-\frac13 a_{22}^2 =a_{51}-\frac23a_{32}a_{22}+\frac13a_{13}a_{22}^2=0.
\end{eqnarray*}
Namely, the last equation $a_{51}-\cdots =0$ is added 
to the condition for  $4_5$-type 
just described above. 
Then, by calculation, one can get $\codim \, G_{W} \ge 4$. 
It is similar for other cases.  
\qed

\

\noindent
{\bf Acknowledgements}
 The author is very grateful to his advisor Toru Ohmoto for a lot of instructions and encouragements. 
Also he thanks K. Saji, T. Nishimura, J. Damon, R. Wik Atique, R. Oset Shinha and M. A. S. Ruas for their comments 
to his talks in workshops and schools held in Hirosaki, Hanoi, S$\tilde{a}$o Calros, 2013-2014.



\begin{thebibliography}{999999}
\bibitem{Arnold}V. I. Arnold, Indices of singular points of 1-forms on  manifolds with boundary, convolution of invariants of groups generated by reflections, and singular projection of smooth hyper surface. Russian Math. Surveys 34 no.2 (1979), 1-42.

%

\bibitem{Arnold2}V. I. Arnold, Singularities of caustics and wavefronts, Kluwer Acad. Publ. (1991).


%



%
\bibitem{Bruce}J. W. Bruce, Projections and reflections of generic surfaces in $\R^3$. Math. Scand. 54 no.2 (1984), 262-278.
%
\bibitem{Damon}J. Damon, Topological triviality and versality for subgroups of $\A$ and $\K$: $I$. Finite codimension conditions. Mem. Amer. Math. Soc. 389 (1988)  
%
\bibitem{Gaffney}
 T. Gaffney, The structure of $T\A(f)$, classification and an application to differential
geometry, In singularities, Part I, Proc. Sympos. in Pure Math. 40 (1983), Amer. Math. Soc., 409-427.
%
\bibitem{Gaffney-Ruas}
 T. Gaffney and M. Ruas, Projections to planes of a geometrically immersed surface (in preparation).
%
\bibitem{Goryunov}
 V. V. Goryunov, Singularities of projections of full intersections, J. Soviet Math. 27. (1984), 
 2785-2811[Translated from Itogi Nauki Tekhniki. Ser. Sovre. Probl. Mat. 22 (1983), 167-206 in Russian].
 %
 \bibitem{Levine}H. I. Levine, The singularities, $S^q_1$, Illinois J. Math. 8 (1964), 152-168.
 %
 
  
  
  
 
 %
%
%


%
\bibitem{ORW1} R. Oset Shinha, M. A. S. Ruas and R. Wik Atique. Classifying codimension two multigerms. Toappear in Mathematical Zeitschrift. Available online DOI: 10.1007/s00209-014-1326-2.

%
%
\bibitem{Platonova}O. A. Platonova, Projections of smooth surfaces, J. Soviet Math.
 35 no. 6 (1986), 2796-2808 [Tr. Sem. I. G. Petvoskii 10 (1984), 135-149 in Russian].

%
\bibitem{Porteous}I. R. Porteous, 
Geometric Differentiation: For the Intelligence of Curves and Surfaces, 2nd edition, 
Cambridge University Press (2001). 
%
\bibitem{Rieger} J. H. Rieger, Families of maps from the plane to the plane.  J. London Math. 
Soc. (2) 36 (1987), no. 2, 351-369.
%
\bibitem{Rieger2} J. H. Rieger, Versal topological stratification and the bifurcation geometry of map-germs of the plane.  Math. Proc. Cambridge Philos.  Soc. 107, no. 1, (1990), 127-147.
%
\bibitem{Rieger3} J. H. Rieger, The geometry of view space of opaque objects bounded by smooth surfaces, Artificial Intelligence. 44 (1990), 1-40.
%

\bibitem{RR}  J. H. Rieger and M. A. S. Ruas, Classification of $\A$-simple germs from $k^n$ to $k^2$, Compositio Math.  79 no. 1, (1991), 99-108. 
%

%
\bibitem{Saji}K. Saji, Criteria for singularities of smooth maps from the plane into the plane and their applications. Hiroshima Math. J. 40, (2010), 229-239.
%
\bibitem{SUY} K. Saji, M. Umehara, and K. Yamada, $A_k$ Singularities of wave fronts. Math. Proc. Cambridge Philos. Soc. 146 (2009), no. 3, 731-746.

\bibitem{SKDO}
H. Sano, Y. Kabata, J. L. Deolindo Silva and T. Ohmoto, Bifurcation in flat geometry of smooth surfaces in $3$-space (in preparation). 


%
\bibitem{Shcherbak}O. P. Shcherbak, Projectively dual curves and Legendre singularities, Selecta Math. Sovietica 5 (1896), 391-421 [Tr. Tbilissk. Univ. 232-233 no13-14 (1982), 280-336 in Russian].

%
\bibitem{Whitney}H. Whitney, On singularities of mappings of Euclidian Spaces I. Mappings of the plane into the plane, Ann. of Math. 62, (1955), 374-410.
%
\bibitem{YKO} T. Yoshida, Y. Kabata and T. Ohmoto, 
Bifurcations of plane-to-plane map-germs of corank $2$, 
Quarterly Jour. Math, (2014), doi:10.1093/qmath/hau013. 
%
\bibitem{YKO2} T. Yoshida, Y. Kabata and T. Ohmoto, 
Bifurcations of plane-to-plane map-germs with corank $2$ of parabolic type, 
accepted in RIMS Bessatsu, Kyoto Univ. (2014). 

\end{thebibliography}
\end{document}